\documentclass[12pt]{article}

\usepackage[T1]{fontenc}
\usepackage{lmodern}
\usepackage{amsmath,amssymb,amsfonts,amsthm,mathtools}
\usepackage{enumitem}
\usepackage{xcolor}
\usepackage{graphicx}
\usepackage{float}
\usepackage{booktabs}
\usepackage{microtype}
\usepackage{tikz}
\usepackage[labelsep=period]{caption}

\usepackage[pdftitle={The strong fractional choice number can exceed the choice number by an arbitrarily large amount},
  pdfauthor={Yangyan Gu and Rongxing Xu},
  pdfkeywords={multiple list coloring, choosability, strong fractional choice number, planar graphs},
  pdfstartview=XYZ,bookmarks=true,colorlinks=true,
  linkcolor=blue,urlcolor=blue,citecolor=blue,
  linktocpage=true,hyperindex=true]{hyperref}

\allowdisplaybreaks
\numberwithin{equation}{section}
\setlist[enumerate]{label=\textup{(\arabic*)},ref=(\arabic*),
  leftmargin=2.35em,itemsep=0.2em,topsep=0.3em}

\newtheorem{theorem}{Theorem}[section]
\newtheorem{lemma}[theorem]{Lemma}

\newtheorem{corollary}[theorem]{Corollary}
\newtheorem{question}[theorem]{Question}
\newtheorem{claim}{Claim}[section]
\newtheorem{observation}{Observation}[section]
\theoremstyle{definition}

\newtheorem{remark}[theorem]{Remark}
\DeclareMathOperator{\ch}{ch}
\newcommand{\chsf}{\ch_f^s}
\newcommand{\N}{\mathbb N}
\newcommand{\set}[1]{\{#1\}}
\newsavebox{\corelistsbox}

\title{\large{\bfseries The strong fractional choice number can exceed the choice number
by an arbitrarily large amount}}
\author{Yangyan Gu\thanks{Shaoxing University, Shaoxing, Zhejiang, 312000, China.
Email: \href{mailto:yygu@usx.edu.cn}{\textcolor{black}{\texttt{yygu@usx.edu.cn}}}.
Supported by the National Natural Science Foundation of China
(Young Scientists Fund, Grant No.~12601684;
General Program, Grant No.~12671412).}
\qquad Rongxing Xu\thanks{School of Mathematical Sciences, Zhejiang Normal University,
Jinhua, Zhejiang, 321004, China.
Email: \href{mailto:xurongxing@zjnu.edu.cn}{\textcolor{black}{\texttt{xurongxing@zjnu.edu.cn}}}.
Supported by the National Natural Science Foundation of China
(Young Scientists Fund, Grant No.~12401472).}
}
\date{\today}

\begin{document}
\maketitle
\vspace{-2em}

\begin{abstract}
Erd\H{o}s, Rubin and Taylor asked whether every $(a,b)$-choosable
graph is $(am,bm)$-choosable for all positive integers $a,b,m$.
Dvo\v{r}\'ak, Hu and Sereni disproved this by constructing a
$4$-choosable graph that is not $(8,2)$-choosable, and asked
whether a $3$-choosable graph can fail to be $(6,2)$-choosable.
Following an earlier discussion by Jiang and Zhu, Li and Zhu asked the
weaker question of whether every $k$-choosable graph is $(kb+1,b)$-choosable
for all positive integers $k,b$.
Equivalently, Xu and Zhu asked whether $\chsf(G)\le\ch(G)$
for every graph $G$. Here $\ch(G)$ is the choice number of a graph $G$, and
$\chsf(G)$ is the strong fractional choice number of $G$,
which is defined as the infimum of the real numbers $r$ such
that $G$ is $(a,b)$-choosable for all positive integers $a,b$
with $a/b\ge r$.

In this paper, we first prove that, for all integers $t\ge0$ and $b\ge7t+2$,
there is a planar graph $G$ with choice number three that is not
$(3b+t,b)$-choosable. For $(b,t)=(2,0)$, we obtain
an affirmative answer to the question of Dvo\v{r}\'ak, Hu and Sereni.
We also answer the questions of Jiang and Zhu, Li and Zhu, and Xu and Zhu
negatively even for planar graphs. The supremum of $\chsf(G)-\ch(G)$
over planar graphs of choice number three is at least $1/7$.
Next, for all integers $k\ge4$, $t\ge0$ and $b\ge3t+2$, we construct a graph $G$ with $\ch(G)=k$ that is not
$(kb+t,b)$-choosable. The case $t=0$ and $b=2$ recovers the main result of
Dvo\v{r}\'ak, Hu and Sereni, and our construction is surprisingly simple,
using only the six-vertex wheel $W_5$ as its core gadget.
The supremum of $\chsf(G)-\ch(G)$ over graphs of choice number $k$
is at least $1/3$ for every integer $k\ge4$.
Finally, for every sufficiently large integer $k$, we construct a graph
$G$ with $\ch(G)=k$ that is not $(2k+t,2)$-choosable, where
$t\ge C\sqrt{k}$ for an absolute constant $C>0$.
Thus the difference between the strong fractional choice number and the
choice number is unbounded even with multiplicity fixed at two.
We further obtain
$k+\Omega(\sqrt{k})\le\sup\set{\chsf(G):\ch(G)=k}\le2^{(1+o(1))k}$
as $k$ tends to infinity.
\end{abstract}

\section{Introduction}\label{sec:introduction}

A \emph{list assignment} of a graph $G$ assigns a set $L(v)$ of colors to each vertex $v\in V(G)$. For a positive integer $b$, a \emph{$b$-fold $L$-coloring} is a mapping $\varphi$ such that $\varphi(v)$ is a $b$-element subset of $L(v)$ for every vertex $v$, and $\varphi(u)\cap\varphi(v)=\varnothing$ for every edge $uv\in E(G)$. For $b=1$, we call this an $L$-coloring. For positive integers $a$ and $b$, a graph is \emph{$(a,b)$-choosable} if it has a $b$-fold $L$-coloring for every list assignment $L$ with $|L(v)|=a$ at every vertex. We abbreviate $(a,1)$-choosable to \emph{$a$-choosable}. The \emph{choice number} $\ch(G)$ is the least integer $a$ such that $G$ is $a$-choosable.

Erd\H{o}s, Rubin and Taylor~\cite{ERT1980} asked whether $(a,b)$-choosability implies $(am,bm)$-choosability for all positive integers $a,b,m$. An affirmative answer is usually called the \emph{$(am,bm)$-conjecture}. Tuza and Voigt~\cite{TuzaVoigt1996} proved that every $(2,1)$-choosable graph is $(2m,m)$-choosable for every positive integer $m$. Dvo\v{r}\'ak, Hu and Sereni~\cite{DvorakHuSereni2019} disproved the conjecture by showing that, for every integer $k\ge4$, there is a $(k,1)$-choosable graph that is not $(2k,2)$-choosable. They posed the following question for $(a,b)=(3,1)$, expecting an affirmative answer.

\begin{question}[Dvo\v{r}\'ak, Hu and Sereni~\cite{DvorakHuSereni2019}]
\label{que:dhs}
Does there exist a graph that is $(3,1)$-choosable but not $(6,2)$-choosable?
\end{question}

Since the original question of Erd\H{o}s, Rubin and Taylor has a negative answer, one may naturally consider weaker versions. Starting from $k$-choosability, the following modification was discussed by Jiang and Zhu~\cite{JiangZhu2019} and posed as a question by Li and Zhu~\cite{LiZhu2020}.

\begin{question}[Jiang and Zhu~\cite{JiangZhu2019}; Li and Zhu~\cite{LiZhu2020}]\label{que:extra-color}
For all positive integers $k,b$, is every $k$-choosable graph $(kb+1,b)$-choosable?
\end{question}

This weakening can be formulated using the strong fractional choice number, introduced by Zhu~\cite{Zhu2017} in 2017. For a real number $r$, a graph $G$ is \emph{strongly fractional $r$-choosable} if it is $(a,b)$-choosable for all positive integers $a,b$ with $a/b\ge r$. The \emph{strong fractional choice number} of $G$ is
\[
 \chsf(G)=\inf\{r\in\mathbb R:
 G\text{ is strongly fractional }r\text{-choosable}\}.
\]
Denote by $\ch_b(G)$ the least integer $a$ such that $G$ is $(a,b)$-choosable. Xu and Zhu~\cite{XuZhu2022} proved that
\begin{equation}\label{eq:strong-choice}
 \chsf(G)=\sup_{b\ge1}\frac{\ch_b(G)-1}{b}.
\end{equation}
In particular, if $G$ is not $(a,b)$-choosable, then $\chsf(G)\ge a/b$. By equation~\eqref{eq:strong-choice}, Question~\ref{que:extra-color} is equivalent to the following question.

\begin{question}[Xu and Zhu~\cite{XuZhu2022}]\label{que:strong-choice}
Does every graph $G$ satisfy $\chsf(G)\le\ch(G)$?
\end{question}

We first answer Questions~\ref{que:dhs}, \ref{que:extra-color}
and~\ref{que:strong-choice}.\footnote{When we
posted the first version of this paper~\cite{GuXu2026v1}, which contained
this result, Hu and Li~\cite{HuLi2026} independently obtained similar
results. Their examples are slightly smaller than ours.
Nevertheless, we obtain planar counterexamples and a
stronger lower bound on the supremum of the strong fractional choice
number over graphs of choice number three.}

\begin{theorem}
\label{prop:planar}
For all integers $t\ge0$ and $b\ge7t+2$, there is a planar graph $G$ with choice number three that is not $(3b+t,b)$-choosable.
As a corollary,
\[
\sup\set{\chsf(G):G\text{ is planar},\ \ch(G)=3}\ge\frac{22}{7}.
\]
\end{theorem}

For $(b,t)=(2,0)$, we obtain a planar $3$-choosable graph that is
not $(6,2)$-choosable, answering Question~\ref{que:dhs}
affirmatively.
For $(b,t)=(9,1)$, we obtain a $3$-choosable graph that is
not $(28,9)$-choosable, so $\chsf(G)\ge28/9>3$.
This answers the questions of Jiang and Zhu, Li and Zhu, and Xu and Zhu
(Questions~\ref{que:extra-color} and~\ref{que:strong-choice}) negatively
even for planar graphs.

The bound $22/7$ follows from $\chsf(G)\ge3+t/b$ by taking $b=7t+2$
and letting $t$ tend to infinity.

For choice number at least four, we prove the following.

\begin{theorem}\label{thm:main}
For all integers $k\ge4$, $t\ge0$ and $b\ge3t+2$, there is a graph $G$ with $\ch(G)=k$ that is not $(kb+t,b)$-choosable.
\end{theorem}

Taking $t=0$ and $b=2$ recovers the result of
Dvo\v{r}\'ak, Hu and Sereni~\cite{DvorakHuSereni2019} that,
for every integer $k\ge4$, there is a $k$-choosable graph that
is not $(2k,2)$-choosable. Our construction is particularly simple,
using only the six-vertex wheel $W_5$ as its core gadget.
For $(k,b,t)=(4,2,0)$, the resulting graph has $2522$ vertices,
compared with $93\,244$ in their construction as originally stated.

Taking $b=3t+2$ in Theorem~\ref{thm:main}, we obtain
$\chsf(G)\ge k+t/(3t+2)$. Letting $t$ tend to infinity, we obtain
the following corollary.

\begin{corollary}\label{cor:third-gap}
For every integer $k\ge4$,
\[
\sup\set{\chsf(G):\ch(G)=k}\ge k+\frac13.
\]
\end{corollary}

The preceding constructions require the multiplicity to grow with the
additive term in the list size. We now fix the multiplicity at two and
obtain a gap of order $\sqrt{k}$ for every sufficiently large choice
number $k$.

\begin{theorem}\label{thm:unbounded-gap}
There is an absolute constant $C>0$ such that, for every sufficiently
large integer $k$, there are a graph $G$ and a positive integer $t$ with
$\ch(G)=k$ and $t\ge C\sqrt{k}$ such that $G$ is not
$(2k+t,2)$-choosable.
\end{theorem}

We obtain the following bounds on the strong fractional choice number
in terms of the choice number.

\begin{corollary}\label{prop:choice-bounds}
As $k$ tends to infinity,
\[
k+\Omega(\sqrt{k})\le\sup\set{\chsf(G):\ch(G)=k}\le2^{(1+o(1))k}.
\]
\end{corollary}

\begin{proof}
For the lower bound, the graph $G$ in Theorem~\ref{thm:unbounded-gap}
satisfies $\chsf(G)\ge k+t/2\ge k+(C/2)\sqrt{k}$, since $G$ is not
$(2k+t,2)$-choosable.

For the upper bound, let $\varepsilon>0$ be fixed.
By Saxton and Thomason~\cite{SaxtonThomason2015}, for all sufficiently
large $k$, every graph of average degree at least $2^{(1+\varepsilon)k}$
has choice number greater than $k$.
Let $D=\lceil2^{(1+\varepsilon)k}\rceil$, and let $G$ be a $k$-choosable graph.
Every nonempty subgraph of $G$ is $k$-choosable and therefore has average
degree less than $D$, so it has a vertex of degree at most $D-1$.
For each positive integer $b$, greedy $b$-fold coloring in the reverse
of a deletion order works from lists of size $Db$, since at most
$(D-1)b$ colors are forbidden at each vertex.
Thus $\chsf(G)\le D$. As $\varepsilon>0$ is arbitrary, we obtain the upper bound.
\end{proof}

In Section~\ref{sec:k3}, we prove Theorem~\ref{prop:planar} using an
eight-vertex graph and record further properties of the planar examples.
In Section~\ref{sec:k4}, we prove Theorem~\ref{thm:main} using the
six-vertex wheel. In Section~\ref{sec:unbounded}, we prove
Theorem~\ref{thm:unbounded-gap} by establishing two coloring properties
of a gadget and then attaching copies of the gadget to a complete graph.

We conclude the introduction with some notation and terminology.
All graphs are finite, simple and nonempty.
For a positive integer $k$, let $[k]=\set{1,2,\ldots,k}$.
For sets $A$ and $B$, let $A\triangle B=(A\setminus B)\cup(B\setminus A)$
denote their \emph{symmetric difference}.
For a vertex $v$ of a graph $G$, let $N_G(v)$ denote its set of neighbors,
and for $S\subseteq V(G)$, write $N_G(S)=\bigcup_{v\in S}N_G(v)$.
Note that $N_G(S)$ may intersect $S$.
We write $\chi(G)$ for the chromatic number of $G$, $\omega(G)$ for its
maximum clique size, and $\alpha(G)$ for its maximum independent set size.

For a function $f\colon V(G)\to\N$, an \emph{$f$-list assignment} has $|L(v)|=f(v)$ at every vertex. The graph $G$ is \emph{$f$-choosable} if it has an $L$-coloring for every $f$-list assignment $L$.

\section{\texorpdfstring{Proof of Theorem~\ref{prop:planar}}{Proof of Theorem 1.4}}\label{sec:k3}\label{sec:core}

Let $H$ be the graph depicted in Figure~\ref{fig:bad-lists}.
Let $f_H\colon V(H)\to\N$ be given by
\[
f_H(w)=\begin{cases}
2,&\text{if }w\in\set{x,y},\\
3,&\text{otherwise}.
\end{cases}
\]

\begin{lemma}\label{lem:ordinary}
The graph $H$ is $f_H$-choosable.
\end{lemma}

\begin{proof}
Let $L$ be an $f_H$-list assignment. It suffices to show that $H$ has an $L$-coloring.

Suppose first that $L(x)\cap L(y)\ne\varnothing$. Color $x,y$ with a common color. At least two colors remain available at every other vertex. Since $H-\set{x,y}$ is a forest, the coloring extends to $H$.

We may therefore assume that $L(x)\cap L(y)=\varnothing$. Since $|L(x)\cup L(y)|=4>|L(u_1)|$, choose $\alpha\in(L(x)\cup L(y))\setminus L(u_1)$. By symmetry, we may assume that $\alpha\in L(x)$. Color $x$ with $\alpha$.

Let $L(y)=\set{\beta_1,\beta_2}$. Suppose that $|L(v_1)\setminus\set{\alpha,\beta_i}|\ge2$ for some $i\in\set{1,2}$. Color $y$ with $\beta_i$. The remaining lists at $v_4,v_3,v_2,v_1$ have sizes at least $1,2,2,2$, respectively, so we can color these vertices greedily in this order and then color $u_2,u_1$. The same argument, with the path colored in reverse order, applies if $|L(v_4)\setminus\set{\alpha,\beta_i}|\ge2$ for some $i\in\set{1,2}$. Since $|L(v_1)|=|L(v_4)|=3$, we may therefore assume that $L(v_1)=L(v_4)=\set{\alpha,\beta_1,\beta_2}$.

Choose $\gamma\in L(v_3)\setminus\set{\beta_1,\beta_2}$, then $c\in L(v_2)\setminus\set{\alpha,\gamma}$. Interchanging $\beta_1$ and $\beta_2$ if necessary, we may assume that $c\ne\beta_1$. Now we color $y$ with $\beta_2$, $v_1,v_4$ with $\beta_1$, $v_2$ with $c$, and $v_3$ with $\gamma$. These colors form a proper partial $L$-coloring of $H$. Since $\alpha\notin L(u_1)$, we can color $u_2$ and then $u_1$ to complete an $L$-coloring of $H$.
\end{proof}

\begin{figure}[h]
\centering
\sbox{\corelistsbox}{%
$\begin{aligned}
 &L(x)=A\cup E,\\
 &L(y)=C\cup F,\\
 &L(u_1)=L(u_2)=E\cup F,\\
 &L(v_1)=A\cup C\cup D\cup E,\\
 &L(v_2)=A\cup B\cup D\cup E,\\
 &L(v_3)=B\cup C\cup D\cup F,\\
 &L(v_4)=A\cup C\cup D\cup F.
\end{aligned}$%
}
\begin{minipage}[t]{.47\linewidth}
\vspace{0pt}
\centering
\resizebox*{!}{\dimexpr\ht\corelistsbox+\dp\corelistsbox\relax}{%
\begin{tikzpicture}[x=.7cm,y=.7cm,
  vertex/.style={circle,fill=black,inner sep=1.5pt},
  every label/.style={font=\small}]
\coordinate (X) at (0,2.4);
\coordinate (Y) at (0,-2.4);
\coordinate (u1) at (-4,0);
\coordinate (u2) at (-2.5,0);
\coordinate (v1) at (-.25,0);
\coordinate (v2) at (1.25,0);
\coordinate (v3) at (2.75,0);
\coordinate (v4) at (4.25,0);
\draw (u1)--(u2);
\draw (v1)--(v2)--(v3)--(v4);
\foreach \a in {u1,u2,v1,v2,v4} \draw (X)--(\a);
\foreach \a in {u1,u2,v1,v3,v4} \draw (Y)--(\a);
\node[vertex,label=above:$x$] at (X) {};
\node[vertex,label=below:$y$] at (Y) {};
\node[vertex,label=below left:$u_1$] at (u1) {};
\node[vertex,label=below:$u_2$] at (u2) {};
\node[vertex,label=below left:$v_1$] at (v1) {};
\node[vertex,label=below:$v_2$] at (v2) {};
\node[vertex,label=above:$v_3$] at (v3) {};
\node[vertex,label=below right:$v_4$] at (v4) {};
\end{tikzpicture}%
}
\end{minipage}\hfill
\begin{minipage}[t]{.50\linewidth}
\vspace{0pt}
\raggedright
\usebox{\corelistsbox}
\end{minipage}
\caption{The graph $H$ and the lists used in the proof of
Lemma~\ref{lem:bad}. The color sets $A,B,C,D,E,F$ are pairwise disjoint.}
\label{fig:bad-lists}
\end{figure}

\begin{lemma}\label{lem:bad}
For all integers $t\ge0$ and $b\ge7t+2$, there is a list assignment $L$ on $H$ such that $|L(w)|=bf_H(w)+t$ for every $w\in V(H)$, and $H$ admits no $b$-fold $L$-coloring.
\end{lemma}

\begin{proof}
Take pairwise disjoint color sets $A,B,C,D,E,F$ with
\[
 |A|=|B|=4t+1,\; |C|=b-3t-1,\; |D|=b-4t-1,\; |E|=2b-3t-1,\; |F|=b+4t+1.
\]
Define the lists as shown in Figure~\ref{fig:bad-lists}.

The lists at $x,y$ have size $2b+t$, and those at $v_2,v_3,v_4,u_1,u_2$ have size $3b+t$. Moreover,
\[
|L(v_1)|=4b-6t-2=(3b+t)+(b-7t-2).
\]
Since $b\ge7t+2$, delete $b-7t-2$ arbitrary colors from $L(v_1)$ and leave all other lists unchanged. The resulting assignment $L$ has the required sizes.

Suppose that $H$ has a $b$-fold $L$-coloring. Count each occurrence of a color at $x$ or $y$ twice, and each occurrence elsewhere once. Since every vertex receives $b$ colors, the total count is $10b$. Deleting colors from $L(v_1)$ can only restrict where a color is used, so it is enough to count using the lists before this deletion.

Each color in $B$ contributes at most one, since $v_2v_3$ is an edge. Each color in $D$ contributes at most two, since $v_1v_2v_3v_4$ is a path. Now consider a color in $A,C,E$ or $F$. In the lists before the deletion, exactly one of $x,y$ has this color in its list, and that vertex is adjacent to every other vertex whose list contains the color. If the color is used there, its contribution is two. Otherwise, after excluding that vertex, the vertices whose lists contain the color induce an edge and an isolated vertex, or two disjoint edges. Its contribution is again at most two. Consequently,
\[
10b\le 2|A|+|B|+2|C|+2|D|+2|E|+2|F|
=10b-1,
\]
a contradiction.
\end{proof}

\begin{proof}[Proof of Theorem~\ref{prop:planar}]
Fix integers $t\ge0$ and $b\ge7t+2$. Let $G$ be obtained from $\binom{3b+t}{b}$ disjoint copies of $H$ by adding a vertex $z$ adjacent to the vertices $x,y$ in each copy.

Given arbitrary lists of size three on $G$, first color $z$. Each copy of $H$ then retains lists of size at least $f_H$. By Lemma~\ref{lem:ordinary}, the coloring extends to every copy, so $\ch(G)\le3$. Since $H$ contains a triangle, $\ch(G)=3$.

Let $L^*$ be the list assignment on $H$ given by Lemma~\ref{lem:bad}. Define a list assignment $L$ of $G$ as follows. Let $L(z)=Z$, where $Z$ is a set of $3b+t$ colors disjoint from all lists of $L^*$. For each $b$-subset $S$ of $Z$, associate a distinct copy of $H$ and assign the lists $L^*(x)\cup S$ and $L^*(y)\cup S$ to its vertices $x$ and $y$, respectively. At every other vertex of this copy, use its list in $L^*$. Thus every list of $L$ has size $3b+t$.

Suppose that $G$ has a $b$-fold $L$-coloring, and let $S$ be the $b$-subset of $Z$ used at $z$. In the copy associated with $S$, no color of $S$ can be used at $x$ or $y$. The restriction of this coloring to the copy is therefore a $b$-fold $L^*$-coloring of $H$, contradicting Lemma~\ref{lem:bad}.

In the embedding of $H$ shown in Figure~\ref{fig:bad-lists}, the vertices $x,y$ lie on a common face. We can add the path $xzy$ in this face and identify all the added vertices $z$ while preserving planarity. Hence $G$ is planar.
\end{proof}

For $(b,t)=(2,0)$, the resulting graph is planar, has $121$ vertices and $240$ edges, and has choice number three while failing to be $(6,2)$-choosable.

\begin{remark}\label{rem:variants}
We can also obtain $2$-connected examples. Use $L^*$ and $Z$ as in the proof of Theorem~\ref{prop:planar}. Take an edge $z_1z_2$ and $\binom{3b+t}{b}\binom{2b+t}{b}$ disjoint copies of $H$, joining $z_1$ to $x$ and $z_2$ to $y$ in each copy. Give both $z_1,z_2$ the list $Z$. For each ordered pair $(S_1,S_2)$ of disjoint $b$-subsets of $Z$, associate a distinct copy of $H$, use its lists $L^*$, and add $S_1$ and $S_2$ to the lists at $x$ and $y$, respectively. For ordinary lists of size three, first color $z_1,z_2$ and then extend the coloring to every copy by Lemma~\ref{lem:ordinary}. In a $b$-fold coloring, the colors chosen at $z_1,z_2$ form one of the indexed pairs $(S_1,S_2)$, so the corresponding copy would have a $b$-fold $L^*$-coloring, contrary to Lemma~\ref{lem:bad}. Since $H$ is $2$-connected and $x,y$ lie on a common face, the resulting graph is $2$-connected and planar. Thus the supremum of $\chsf(G)$ over $2$-connected planar graphs of choice number three is also at least $22/7$.
\end{remark}

\section{\texorpdfstring{Proof of Theorem~\ref{thm:main}}{Proof of Theorem 1.5}}\label{sec:k4}

Let $W$ be the wheel graph depicted in Fig.~\ref{fig:wheel}, with center $u$ and rim $v_0v_1v_2v_3v_4v_0$. Define $f_W\colon V(W)\to\N$ by
\[
f_W(u)=4,\quad f_W(v_0)=f_W(v_1)=f_W(v_2)=2,\quad f_W(v_3)=f_W(v_4)=3.
\]

\begin{figure}[ht]
\centering
\begin{tikzpicture}[scale=.7,vertex/.style={circle,fill=black,inner sep=1.35pt},every label/.style={font=\small}]
\coordinate (u) at (0,0);
\foreach \i in {0,...,4}
  \coordinate (v\i) at ({162-72*\i}:2);
\foreach \a/\b in {v0/v1,v1/v2,v2/v3,v3/v4,v4/v0,u/v0,u/v1,u/v2,u/v3,u/v4} \draw (\a)--(\b);
\node[vertex,label=below:$u$] at (u) {};
\node[vertex,label=above left:$v_0$] at (v0) {};
\node[vertex,label=above:$v_1$] at (v1) {};
\node[vertex,label=above right:$v_2$] at (v2) {};
\node[vertex,label=below right:$v_3$] at (v3) {};
\node[vertex,label=below left:$v_4$] at (v4) {};
\end{tikzpicture}
\caption{The wheel $W_5$.}
\label{fig:wheel}
\end{figure}

\begin{lemma}\label{lem:wheel-ordinary}
The graph $W$ is $f_W$-choosable.
\end{lemma}
\begin{proof}
Let $L$ be an $f_W$-list assignment. It suffices to show that $W$ has an $L$-coloring.

Suppose first that $L(v_3)\ne L(v_4)$. Choose $\alpha\in L(v_4)\setminus L(v_3)$. We can greedily color $v_0,v_1,v_2,u$, in this order, to obtain a proper partial $L$-coloring $\phi$ with $\phi(v_0)\ne\alpha$. If $\phi(u)=\alpha$, at least two colors remain available at $v_3$, since $\alpha\notin L(v_3)$, and at least one color remains available at $v_4$. Thus we can color $v_4$ first and then $v_3$. If $\phi(u)\ne\alpha$, set $\phi(v_4)=\alpha$, which differs from both $\phi(u)$ and $\phi(v_0)$. Since $\alpha\notin L(v_3)$, at least one color remains available at $v_3$, completing the coloring.

We may therefore assume that $L(v_3)=L(v_4)$. Suppose now that $L(v_0)\ne L(v_1)$. Choose $\alpha_0\in L(v_0)\setminus L(v_1)$, then $\alpha_2\in L(v_2)\setminus\set{\alpha_0}$, and finally $\alpha_1\in L(v_1)\setminus\set{\alpha_2}$. Note that $\alpha_0,\alpha_1,\alpha_2$ are pairwise distinct. Now we color $v_i$ with $\alpha_i$ for $i=0,1,2$ and $u$ with a color $\beta\in L(u)\setminus\set{\alpha_0,\alpha_1,\alpha_2}$. The remaining lists at $v_3,v_4$ are $L(v_3)\setminus\set{\beta,\alpha_2}$ and $L(v_4)\setminus\set{\beta,\alpha_0}$, respectively. Both remaining lists are nonempty, since each original list has size three. If $\max\set{|L(v_3)\setminus\set{\beta,\alpha_2}|,|L(v_4)\setminus\set{\beta,\alpha_0}|}\ge2$, color a vertex with a smallest remaining list first and then the other vertex. Otherwise, $|L(v_3)\setminus\set{\beta,\alpha_2}|=|L(v_4)\setminus\set{\beta,\alpha_0}|=1$, so $\set{\beta,\alpha_2}\subseteq L(v_3)$ and $\set{\beta,\alpha_0}\subseteq L(v_4)$. Since $L(v_3)=L(v_4)$ has size three, this common list is $\set{\beta,\alpha_0,\alpha_2}$. Now we color $v_3$ with $\alpha_0$ and $v_4$ with $\alpha_2$. These colors are available and distinct, so this completes a proper $L$-coloring of $W$.

By symmetry, $W$ is also $L$-colorable if $L(v_2)\ne L(v_1)$. We may therefore assume that $L(v_0)=L(v_1)=L(v_2)=\set{\alpha,\beta}$. If $\set{\alpha,\beta}\subseteq L(v_3)$, color $u$ with a color $\gamma\in L(u)\setminus L(v_3)$, color $v_0,v_2$ with $\alpha$, and color $v_1$ with $\beta$. Then two colors remain available at each of $v_3,v_4$, so we can color $v_3$ and then $v_4$. Otherwise, we may assume that $\alpha\notin L(v_3)$. Color $v_0,v_2$ with $\alpha$, color $v_1$ with $\beta$, and color $u$ with a color $\gamma\in L(u)\setminus\set{\alpha,\beta}$. At least two colors remain available at each of $v_3,v_4$, so again we can color $v_3$ and then $v_4$.
\end{proof}

\begin{lemma}\label{lem:wheel-bad}
For all integers $t\ge0$ and $b\ge3t+2$, there exists a list assignment $L^*$ on $W$ such that $|L^*(x)|=bf_W(x)+t$ for every $x\in V(W)$, and $W$ admits no $b$-fold $L^*$-coloring.
\end{lemma}
\begin{proof}
Take pairwise disjoint color sets $X_0,\ldots,X_4,Y$ such that
\[
|X_0|=|X_1|=|X_2|=|X_4|=b-2t-1,\; |X_3|=2b-2t-1,\; |Y|=5t+2.
\]
For $i=0,\ldots,4$, let
\[
L^*(v_i)=X_{i-1}\cup X_i\cup Y,
\]
where the subscripts of $v_i$ and $X_i$ are taken modulo $5$. Since $b\ge3t+2$, the union of these six color sets has size
\[
6b-5t-3=4b+t+2(b-3t-2)+1>4b+t.
\]
Choose $L^*(u)$ to be any $(4b+t)$-element subset of $X_0\cup\cdots\cup X_4\cup Y$. The resulting lists have the required sizes.

Suppose that $W$ has a $b$-fold $L^*$-coloring. Each color in $X_i$ can be used only at vertices of the triangle $uv_iv_{i+1}$, so it contributes at most one occurrence. Each color in $Y$ contributes at most two occurrences, since every independent set of $W$ has size at most two. Consequently,
\[
6b\le\sum_{i=0}^4|X_i|+2|Y|=6b-1,
\]
a contradiction.
\end{proof}

\begin{proof}[Proof of Theorem~\ref{thm:main}]
Let $k\ge4$, $t\ge0$ and $b\ge3t+2$ be fixed integers.
Take a complete graph $Q$ on vertices $z_1,\ldots,z_{k-2}$
and $N$ pairwise disjoint copies $W^1,\ldots,W^N$ of $W$, all disjoint from $Q$, where
\[
N=\binom{kb+t}{b}\binom{(k-1)b+t}{b}\cdots\binom{3b+t}{b}.
\]
For each $j\in[N]$, the center of $W^j$ is labeled $u^j$ and the rim vertices are labeled $v_0^j,\ldots,v_4^j$ in the cyclic order shown in Fig.~\ref{fig:wheel}.
For every $j\in[N]$, add the following edges between $V(W^j)$ and $V(Q)$:
\begin{enumerate}
\item Join $u^j$ to each of $z_1,\ldots,z_{k-4}$.
\item Join each of $v_3^j,v_4^j$ to each of $z_1,\ldots,z_{k-3}$.
\item Join each of $v_0^j,v_1^j,v_2^j$ to each of $z_1,\ldots,z_{k-2}$.
\end{enumerate}
When $k=4$, each $u^j$ has no neighbor in $Q$.
No other edges are added; in particular, distinct copies of $W$ have no edges between them. Let $G$ be the resulting graph; see Fig.~\ref{fig:construction}.

\begin{figure}[htbp]
\centering
\colorlet{qouterfill}{gray!20}
\colorlet{qmiddlefill}{cyan!20}
\colorlet{qinnerfill}{magenta!20}
\begin{tikzpicture}[
  xscale=1.45,
  vertex/.style={circle,fill=black,inner sep=1.6pt},
  attachment/.style={draw=black,line width=1.05pt},
  redattachment/.style={draw=red,line width=1.05pt},
  blueattachment/.style={draw=blue,line width=1.05pt},
  every label/.style={font=\small,inner sep=1pt},
  line cap=round]
\fill[fill=qouterfill,rounded corners=2pt]
  (-4.6,2.60) rectangle (4.6,4.20);
\fill[fill=qmiddlefill,rounded corners=2pt]
  (-4.45,2.80) rectangle (2.45,4.00);
\fill[fill=qinnerfill,rounded corners=2pt]
  (-4.3,3.00) rectangle (.3,3.80);
\node[anchor=west,inner sep=2pt] at (4.8,3.40) {$Q$};
\coordinate (Zfirst) at (-4,3.25);
\coordinate (Zsecond) at (-2,3.25);
\coordinate (Zlast) at (0,3.25);
\coordinate (Znext) at (2,3.25);
\coordinate (Zend) at (4,3.25);
\coordinate (CU) at (-4.25,0);
\coordinate (CV3) at (-2.55,0);
\coordinate (CV4) at (-.85,0);
\coordinate (CV0) at (.85,0);
\coordinate (CV1) at (2.55,0);
\coordinate (CV2) at (4.25,0);
\draw[draw=black!45,fill=black!1,line width=.5pt,rounded corners=2pt]
  (-4.6,-.65) rectangle (4.6,.28);
\node[anchor=west,inner sep=2pt] at (4.8,-.185) {$W^j$};
\node at (-5.2,-.185) {$\cdots$};
\node at (5.75,-.185) {$\cdots$};
\foreach \z in {Zfirst,Zsecond,Zlast}
  \draw[redattachment] (CU)--(\z);
\foreach \i in {0,1,2}
  \foreach \z in {Zfirst,Zsecond,Zlast,Znext,Zend}
    \draw[attachment] (CV\i)--(\z);
\foreach \i in {3,4}
  \foreach \z in {Zfirst,Zsecond,Zlast,Znext}
    \draw[blueattachment] (CV\i)--(\z);
\node[vertex,label=above:$z_1$] at (Zfirst) {};
\node[vertex,label=above:$z_2$] at (Zsecond) {};
\node at (-1,3.25) {$\cdots$};
\node[vertex,label=above:$z_{k-4}$] at (Zlast) {};
\node[vertex,label=above:$z_{k-3}$] at (Znext) {};
\node[vertex,label=above:$z_{k-2}$] at (Zend) {};
\node[vertex,label=below:$u^j$] at (CU) {};
\foreach \i in {0,...,4}
  \node[vertex,label=below:$v_\i^j$] at (CV\i) {};
\end{tikzpicture}
\caption{The construction of $G$, with the internal edges of $Q$ and $W^j$
omitted.}
\label{fig:construction}
\end{figure}

We first prove that $\ch(G)=k$. Let $L$ be a $k$-list assignment of $G$.
First color $Q$ greedily from $L$, and delete from the list of each vertex outside $Q$ the colors used on its neighbors in $Q$.
Thus, for every $j\in[N]$, at least $k-(k-4)=4$ colors remain at $u^j$,
at least $k-(k-3)=3$ at each of $v_3^j,v_4^j$, and at least $k-(k-2)=2$ at
each of $v_0^j,v_1^j,v_2^j$. By Lemma~\ref{lem:wheel-ordinary}, the
coloring extends to every $W^j$. Since distinct copies have no edges between them, we obtain a proper $L$-coloring of $G$, so $\ch(G)\le k$.
The vertices $z_1,\ldots,z_{k-4}$ form a clique adjacent to
every vertex of each $W^j$. Since $\chi(W)=4$, we have
$\chi(G)\ge(k-4)+4=k$, and hence $\ch(G)=k$.

We next prove that $G$ is not $(kb+t,b)$-choosable. Let $L^*$ be the list assignment on $W$ given by Lemma~\ref{lem:wheel-bad}, and let $Z$ be a set of $kb+t$ colors disjoint from $\bigcup_{x\in V(W)}L^*(x)$. We construct a $(kb+t)$-list assignment $L$ of $G$, beginning with $L(z_i)=Z$ for $i\in[k-2]$.

The $b$-fold $L$-colorings of $Q$ correspond to the ordered tuples $(S_1,\ldots,S_{k-2})$ of pairwise disjoint $b$-subsets of $Z$, with $S_i$ used at $z_i$. There are exactly $N$ such tuples, so we associate them bijectively with the copies $W^1,\ldots,W^N$. In the copy $W^j$ associated with $(S_1,\ldots,S_{k-2})$, we add $S_i$ to the original list of each neighbor of $z_i$ in this copy, as follows:
\begin{itemize}
\item For $u^j$, let $L(u^j)=L^*(u)\cup\bigcup_{i=1}^{k-4}S_i$.
\item For $\ell\in\set{3,4}$, let $L(v_\ell^j)=L^*(v_\ell)\cup\bigcup_{i=1}^{k-3}S_i$.
\item For $\ell\in\set{0,1,2}$, let $L(v_\ell^j)=L^*(v_\ell)\cup\bigcup_{i=1}^{k-2}S_i$.
\end{itemize}
Clearly, $|L(x^j)|=kb+t$ for every $x\in V(W)$.

Suppose that $G$ has a $b$-fold $L$-coloring $\phi$, and let $S_i=\phi(z_i)$ for $i\in[k-2]$. Since $Q$ is complete, these sets form an ordered tuple of pairwise disjoint $b$-subsets of $Z$. Consider the copy $W^j$ associated with this tuple. For each $x\in V(W)$, every set $S_i$ added to $L^*(x)$ in forming $L(x^j)$ is used at a neighbor $z_i$ of $x^j$. None of these added colors can therefore be used at $x^j$, so $\phi(x^j)\subseteq L^*(x)$. The restriction of $\phi$ to $W^j$ is consequently a $b$-fold $L^*$-coloring of $W$, contradicting Lemma~\ref{lem:wheel-bad}.
\end{proof}

\section{\texorpdfstring{Proof of Theorem~\ref{thm:unbounded-gap}}{Proof of Theorem 1.7}}\label{sec:unbounded}

We construct the auxiliary graph $J$ in Subsection~\ref{subsec:twofold-auxiliary}.
In Subsection~\ref{subsec:twofold-construction}, we build $H(t,d)$ from $J$.
For each vertex $v$ of $H(t,d)$, we specify a list size $f(v)$ for ordinary
list coloring. We also construct a list assignment that admits no $2$-fold
coloring and has $2f(v)+t$ colors at each vertex $v$. We prove that
$H(t,d)$ is $f$-choosable when $28000t<2^{d/2-9}$.
We complete the proof in Subsection~\ref{subsec:twofold-uniform} by
attaching copies of $H(t,d)$ to a complete graph.

To prove $f$-choosability, we color pairs of nonadjacent vertices of $J$
with common colors. If we cannot choose distinct colors from the
remaining lists, Hall's theorem identifies a set of vertices whose
lists together contain fewer colors than vertices. Each color then
determines a clique in $J$. We approximate these cliques by cliques indexed
by binary vectors and use the resulting counts to define an integer-valued
function whose signed sums
contradict the parity estimate in Subsection~\ref{subsec:binary-estimates}.

\subsection{Signed sums for prescribed parities}\label{subsec:binary-estimates}

Let $d$ be a positive integer. For
$\mathbf{u}=(u_1,\ldots,u_d)$ and $\mathbf{x}=(x_1,\ldots,x_d)$ in
$\mathbb F_2^d$, we define their dot product by
\[
\mathbf{u}\cdot\mathbf{x}
=\left(\sum_{i=1}^d u_i x_i\right)\bmod 2.
\]
For a real-valued function $g$ on $\mathbb F_2^d$ and
$\mathbf{x}\in\mathbb F_2^d$, let
\[
\Delta_{\mathbf{x}}(g)=
\sum_{\mathbf{u}\in\mathbb F_2^d}g(\mathbf{u})
(-1)^{\mathbf{u}\cdot\mathbf{x}}.
\]
Thus $\Delta_{\mathbf{x}}(g)$ is the sum of $g(\mathbf{u})$ over
the vectors with $\mathbf{u}\cdot\mathbf{x}=0$ minus the sum
over the vectors with $\mathbf{u}\cdot\mathbf{x}=1$.

\begin{lemma}\label{lem:parseval}
Let $g$ be a real-valued function on $\mathbb F_2^d$. Then
\[
\sum_{\mathbf{x}\in\mathbb F_2^d}\Delta_{\mathbf{x}}(g)^2
=2^d\sum_{\mathbf{u}\in\mathbb F_2^d}g(\mathbf{u})^2.
\]
\end{lemma}

\begin{proof}
Expanding the squares and using the additivity of the dot product,
we obtain
\[
\sum_{\mathbf{x}\in\mathbb F_2^d}\Delta_{\mathbf{x}}(g)^2
=\sum_{\mathbf{u},\mathbf{v}\in\mathbb F_2^d}
g(\mathbf{u})g(\mathbf{v})
\sum_{\mathbf{x}\in\mathbb F_2^d}
(-1)^{(\mathbf{u}+\mathbf{v})\cdot\mathbf{x}}.
\]
If $\mathbf{u}=\mathbf{v}$, then $\mathbf{u}+\mathbf{v}=\mathbf{0}$,
so all $2^d$ terms of the inner sum are $1$. 
If $\mathbf{u}\ne\mathbf{v}$, let $i$ be the first coordinate with
$u_i+v_i=1$. We pair the vectors so that the two vectors in each pair
differ only in coordinate $i$. In each pair, let $\mathbf{x}$ be the
vector with $x_i=0$ and $\mathbf{x}'$ the vector with $x_i'=1$.
Since $u_i+v_i=1$, we have
$(\mathbf{u}+\mathbf{v})\cdot\mathbf{x}'
=(\mathbf{u}+\mathbf{v})\cdot\mathbf{x}+1$ in $\mathbb F_2$.
The two terms in each pair cancel, so the inner sum is $0$. 
The right-hand side therefore reduces to
$2^d\sum_{\mathbf{u}\in\mathbb F_2^d}g(\mathbf{u})^2$.
\end{proof}

\begin{lemma}\label{lem:fourier-parity}
Let $d\ge6$ be an integer. There is an even-sized subset $B$ of
$\mathbb F_2^d$ such that every integer-valued function $g$ on
$\mathbb F_2^d$ that is odd exactly on $B$ satisfies
\[
\sum_{\mathbf{x}\in\mathbb F_2^d}|\Delta_{\mathbf{x}}(g)|\ge2^{3d/2-9}.
\]
\end{lemma}

\begin{proof}
Let $\mathcal G$ be the set of integer-valued functions $g$ on
$\mathbb F_2^d$ with
$\sum_{\mathbf{x}\in\mathbb F_2^d}|\Delta_{\mathbf{x}}(g)|<2^{3d/2-9}$, and 
\[
P_g=\{\mathbf{u}\in\mathbb F_2^d\mid g(\mathbf{u})\equiv1\pmod{2}\}
\quad(g\in\mathcal G),
\qquad 
\mathcal F=\{P_g\mid g\in\mathcal G\}.
\]
Let $\tau=2^{d/2-3}$ and $s=2^{d-6}$. For $g\in\mathcal G$ and
$\mathbf{x}\in\mathbb F_2^d$, write $z_{\mathbf{x}}(g)$ for
$\Delta_{\mathbf{x}}(g)/\tau$ rounded toward zero. Then
\[
\sum_{\mathbf{x}\in\mathbb F_2^d}|z_{\mathbf{x}}(g)|
\le\frac{1}{\tau}\sum_{\mathbf{x}\in\mathbb F_2^d}|\Delta_{\mathbf{x}}(g)|
<\frac{2^{3d/2-9}}{2^{d/2-3}}=s.
\]

For each $g \in \mathcal G$, we collect the values $z_{\mathbf{x}}(g)$ into a tuple. Let
$$\mathcal T_1=\bigl\{(z_{\mathbf{x}}(g))_{\mathbf{x}\in\mathbb F_2^d}
\mid g\in\mathcal G\bigr\} \quad \text{and} \quad \mathcal T_2=\bigl\{(|z_{\mathbf{x}}(g)|)_{\mathbf{x}\in\mathbb F_2^d}
\mid g\in\mathcal G\bigr\}.$$

We group the functions in $\mathcal G$ by their tuples in $\mathcal T_1$,
so there are $|\mathcal T_1|$ groups. The entries of a tuple in
$\mathcal T_2$ sum to $\sum_{\mathbf{x}}|z_{\mathbf{x}}(g)|<s$,
so at most $s$ entries are nonzero. Each nonzero entry has two possible
signs, so each tuple in $\mathcal T_2$ comes from at most $2^s$ tuples
in $\mathcal T_1$.

There are $\binom{2^d+j-1}{j}$ nonnegative integer $2^d$-tuples with
sum $j$. By Pascal's identity,
\[
|\mathcal T_2|\le\sum_{j=0}^{s}\binom{2^d+j-1}{j}
=1+\sum_{j=1}^{s}\left(\binom{2^d+j}{j}-\binom{2^d+j-1}{j-1}\right)
=\binom{2^d+s}{s}.
\]
Using $\binom{n}{j}\le(en/j)^j$ and $2^6s=2^d$, we obtain
\begin{equation}\label{eq:parity-group-count}
|\mathcal T_1|\le 2^s|\mathcal T_2|\le 2^s\binom{2^d+s}{s}
\le 2^s\left(\frac{e(2^d+s)}{s}\right)^s
=(130e)^s<2^{16s}.
\end{equation}

In each nonempty group, let $g_0$ be a fixed member. For each $g$ in
the group and each $\mathbf{x}\in\mathbb F_2^d$,
both $\Delta_{\mathbf{x}}(g)$ and $\Delta_{\mathbf{x}}(g_0)$ differ from
$\tau z_{\mathbf{x}}(g)$ by less than $\tau$, so
$|\Delta_{\mathbf{x}}(g)-\Delta_{\mathbf{x}}(g_0)|<2\tau$.
By Lemma~\ref{lem:parseval} applied to $g-g_0$, we have
\[
\sum_{\mathbf{u}\in\mathbb F_2^d}
\bigl(g(\mathbf{u})-g_0(\mathbf{u})\bigr)^2
=\frac{1}{2^d}\sum_{\mathbf{x}\in\mathbb F_2^d}
\bigl(\Delta_{\mathbf{x}}(g)-\Delta_{\mathbf{x}}(g_0)\bigr)^2
<4\tau^2=2^{d-4}.
\]

For each $\mathbf{u}\in P_g\triangle P_{g_0}$, one of $g(\mathbf{u})$
and $g_0(\mathbf{u})$ is odd and the other is even, so the squared
difference is at least $1$.
Thus $|P_g\triangle P_{g_0}|<2^{d-4}=4s$.
Since $(P_g\triangle P_{g_0})\triangle P_{g_0}=P_g$, each possible
symmetric difference determines at most one set $P_g$.
Thus each group accounts for at most
\[
\sum_{j=0}^{4s}\binom{2^d}{j}
\le 16^{4s}\sum_{j=0}^{2^d}\binom{2^d}{j}16^{-j}
=16^{4s}(1+1/16)^{2^d}
\le (16e)^{4s}<2^{24s}
\]
sets $P_g$. Here $16^{4s-j}\ge1$ for $0\le j\le4s$,
$1+1/16\le e^{1/16}$, and $2^d=64s$.

By inequality~\eqref{eq:parity-group-count}, there are fewer than $2^{16s}$ groups,
so $|\mathcal F|<2^{16s}\cdot 2^{24s}=2^{40s}$.

Suppose the conclusion does not hold. Then for each even-sized subset
$B\subseteq\mathbb F_2^d$, there exists $g\in\mathcal G$ with $P_g=B$,
so $B\in\mathcal F$.
Note that the number of even-sized subsets of $\mathbb F_2^d$ is
\[
\sum_{j=0}^{2^{d-1}}\binom{2^d}{2j}
=\frac{(1+1)^{2^d}+(1-1)^{2^d}}{2}
=2^{2^d-1}=2^{64s-1}.
\]
Thus $|\mathcal F|\ge2^{64s-1}>2^{40s}$, a contradiction.
\end{proof}

\subsection{\texorpdfstring{Construction of an auxiliary graph $J$}{Construction of an auxiliary graph J}}\label{subsec:twofold-auxiliary}

Let $d\ge6$ be even. In $\mathbb F_2^d$, let $\Omega_1$ consist of the vectors
whose last $d/2$ coordinates are zero, and let $\Omega_2$ consist of the vectors
whose first $d/2$ coordinates are zero. Each vector is uniquely a sum of
one vector in $\Omega_1$ and one in $\Omega_2$. 
We use lexicographic order on $\mathbb F_2^d$, comparing vectors at the
first coordinate where they differ, with $0<1$.
Let
\[
\begin{aligned}
\mathcal T={}&
\bigl\{(\mathbf{x},\mathbf{y},\mathbf{x}+\mathbf{y})
\mid \mathbf{x}<\mathbf{y}<\mathbf{x}+\mathbf{y},\
\text{$\mathbf{x},\mathbf{y}$ both in $\Omega_1\setminus\{\mathbf{0}\}$
or both in $\Omega_2\setminus\{\mathbf{0}\}$}\bigr\}\\
&\cup\bigl\{(\mathbf{x},\mathbf{y})
\mid \mathbf{x}\in \Omega_1\setminus\{\mathbf{0}\},\
\mathbf{y}\in \Omega_2\setminus\{\mathbf{0}\}\bigr\}.
\end{aligned}
\]
Note that each triple $(\mathbf{x},\mathbf{y},\mathbf{z})$ of vectors all
in $\Omega_1$ or all in $\Omega_2$ with $\mathbf{0}<\mathbf{x}<\mathbf{y}<\mathbf{z}$
and $\mathbf{x}+\mathbf{y}+\mathbf{z}=\mathbf{0}$ determines a unique
member of $\mathcal T$ with these entries. Indeed,
$\mathbf{x}+\mathbf{y}=(\mathbf{x}+\mathbf{y}+\mathbf{z})+\mathbf{z}=\mathbf{z}$.
Thus the triple in $\mathcal T$ is
$(\mathbf{x},\mathbf{y},\mathbf{x}+\mathbf{y})$, with its entries in the
only increasing order.

\begin{observation}\label{obs:tuple-count}
$|\mathcal T|<(4/3)2^d$.
\end{observation}

\begin{proof}
Each of $\Omega_1$ and $\Omega_2$ has $2^{d/2}-1$ nonzero vectors. In either subspace,
each ordered pair of distinct nonzero vectors determines one triple,
and each triple comes from six such pairs. A two-entry tuple has one
nonzero vector from $\Omega_1$ and one from $\Omega_2$. Thus
\[
|\mathcal T|
=2\cdot\frac{(2^{d/2}-1)(2^{d/2}-2)}{6}+(2^{d/2}-1)^2
<\frac43\cdot2^d.\qedhere
\]
\end{proof}

For $T\in\mathcal T$ and $\mathbf{u}\in\mathbb F_2^d$, let
\[
\beta_{\mathbf{u}}(T)=
\begin{cases}
(\mathbf{u}\cdot\mathbf{x},\mathbf{u}\cdot\mathbf{y}),
&\text{if } T=(\mathbf{x},\mathbf{y}),\\
(\mathbf{u}\cdot\mathbf{x},\mathbf{u}\cdot\mathbf{y},\mathbf{u}\cdot(\mathbf{x}+\mathbf{y})),
&\text{if } T=(\mathbf{x},\mathbf{y},\mathbf{x}+\mathbf{y}).
\end{cases}
\]
For each $T\in\mathcal T$, we say that two vectors
$\mathbf{u},\mathbf{v}\in\mathbb F_2^d$ are equivalent with respect to $T$
if $\beta_{\mathbf{u}}(T)=\beta_{\mathbf{v}}(T)$.

\begin{observation}\label{obs:dot-product-partition}
Each $T\in\mathcal T$ has exactly four equivalence classes, each of size $2^{d-2}$.
\end{observation}

\begin{proof}
Let $\mathbf{x}=(x_1,\ldots,x_d)$ and $\mathbf{y}=(y_1,\ldots,y_d)$
be the first two entries of $T$. The first two entries of
$\beta_{\mathbf{u}}(T)$ determine the whole tuple, since for a triple
the third entry is their sum. Thus for each $(i,j)\in\mathbb F_2^2$,
we count the vectors $\mathbf{u}=(u_1,\ldots,u_d)$ satisfying
$\mathbf{u}\cdot\mathbf{x}=i$ and $\mathbf{u}\cdot\mathbf{y}=j$.

If $T$ is a triple, let $p$ be the first index with $x_p\ne y_p$.
Since $\mathbf{x}<\mathbf{y}$ in lexicographic order, $x_p=0$ and $y_p=1$.
If $T$ is a pair, the nonzero vector $\mathbf{y}\in \Omega_2$ has $y_p=1$ at
some index $p>d/2$. At that index $x_p=0$ because $\mathbf{x}\in \Omega_1$.
In either case, there is an index $q\ne p$ with $x_q=1$ because
$\mathbf{x}\ne\mathbf{0}$.

The entries $u_k$ with $k\notin\{p,q\}$ can be chosen freely. Since
$x_p=0$ and $x_q=1$, the equation $\mathbf{u}\cdot\mathbf{x}=i$
determines $u_q$. Then $\mathbf{u}\cdot\mathbf{y}=j$ determines $u_p$
because $y_p=1$. Thus each of the four pairs $(i,j)$ occurs for exactly
$2^{d-2}$ vectors and corresponds to one equivalence class.
\end{proof}

The vertices of $J$ are the pairs $v=(T,X)$, where $T\in\mathcal T$
and $X$ is an equivalence class with respect to $T$. For each such vertex,
let $X_v=X$. Distinct vertices $v,v'$ are adjacent if and only if
$X_v\cap X_{v'}\ne\varnothing$.

\begin{observation}\label{obs:auxiliary-graph-properties}
$|V(J)|=4|\mathcal T|<6\cdot2^d$ and $\alpha(J)=4$.
\end{observation}

\begin{proof}
By Observation~\ref{obs:dot-product-partition}, each $T\in\mathcal T$
corresponds to four vertices of $J$, so $|V(J)|=4|\mathcal T|$.
By Observation~\ref{obs:tuple-count}, $|V(J)|<6\cdot2^d$.

For a fixed $T$, its four equivalence classes are pairwise disjoint,
so the corresponding vertices form an independent set and $\alpha(J)\ge4$.
If $I$ is an independent set in $J$, then $X_v\cap X_{v'}=\varnothing$
for distinct $v,v'\in I$. These sets are contained in $\mathbb F_2^d$
and each has size $2^{d-2}$ by Observation~\ref{obs:dot-product-partition}.
Thus $|I|\cdot 2^{d-2}\le2^d$, so $\alpha(J)\le4$.
\end{proof}

A function $\varphi\colon\mathbb F_2^n\to\mathbb F_2$ is \emph{linear} if
$\varphi(\mathbf{x}+\mathbf{y})=\varphi(\mathbf{x})+\varphi(\mathbf{y})$
for all $\mathbf{x},\mathbf{y}\in\mathbb F_2^n$.
We use the following lemma. For a proof, see O'Donnell~\cite{ODonnell2014}
(p.~15, Theorem~1.30).

\begin{lemma}\label{lem:binary-linearity}
Let $\varepsilon\in[0,1]$ and let $\psi\colon\mathbb F_2^n\to\mathbb F_2$. If
$\psi(\mathbf{x})+\psi(\mathbf{y})\ne\psi(\mathbf{x}+\mathbf{y})$
for at most $\varepsilon 2^{2n}$ choices of
$\mathbf{x},\mathbf{y}\in\mathbb F_2^n$, then there is a linear function
$\varphi\colon\mathbb F_2^n\to\mathbb F_2$ that differs from $\psi$ on at most
$\varepsilon 2^n$ points.
\end{lemma}

For each $\mathbf{u}\in\mathbb F_2^d$, let
$Q_{\mathbf{u}}=\{v\in V(J):\mathbf{u}\in X_v\}$. The vertices of
$Q_{\mathbf{u}}$ form a clique in $J$ because their sets $X_v$ all
contain $\mathbf{u}$.
\begin{lemma}\label{lem:clique-approximation}
	Let $C\subseteq V(J)$ be a clique, possibly empty, and let
	$r=|\mathcal T|-|C|$. There is a vector
	$\mathbf{u}\in\mathbb F_2^d$ such that
	\[
	|C\triangle Q_{\mathbf{u}}|\le19r.
	\]
\end{lemma}

\begin{proof}
	The vertices with the same $T$ are nonadjacent, so $C$ contains at most one vertex $(T,X)$ for each $T$. Thus $r$ counts the tuples $T\in\mathcal T$ for which $C$ contains no such vertex.
	
	For each $(T,X)\in C$ and each entry $\mathbf{x}$ of $T$, define $\psi(\mathbf{x})=\mathbf{w}\cdot\mathbf{x}$ for any $\mathbf{w}\in X$. This value does not depend on the choice of $\mathbf{w}$, by the definition of the equivalence class $X$. If $\mathbf{x}$ also occurs in $T'$ for another vertex $(T',X')\in C$, then $X\cap X'\ne\varnothing$, since $C$ is a clique. For any $\mathbf{w}\in X\cap X'$, both definitions give $\psi(\mathbf{x})=\mathbf{w}\cdot\mathbf{x}$. Thus $\psi$ is well defined.
	For every vector $\mathbf{x}\in\Omega_1\cup\Omega_2$ that occurs in none of the tuples represented in $C$, set $\psi(\mathbf{x})=0$. In particular, $\psi(\mathbf{0})=0$.
	
	For $\mathbf{x},\mathbf{y}$ both in $\Omega_1$ or both in $\Omega_2$, the equality
	\(
	\psi(\mathbf{x})+\psi(\mathbf{y})=\psi(\mathbf{x}+\mathbf{y})
	\)
	holds if either vector is zero or $\mathbf{x}=\mathbf{y}$. For distinct nonzero $\mathbf{x},\mathbf{y}$, let $T$ be the triple with entries $\mathbf{x},\mathbf{y},\mathbf{x}+\mathbf{y}$. If $(T,X)\in C$, then any $\mathbf{w}\in X$ gives
	\[
	\psi(\mathbf{x})+\psi(\mathbf{y})
	=\mathbf{w}\cdot(\mathbf{x}+\mathbf{y})
	=\psi(\mathbf{x}+\mathbf{y}).
	\]
	
	For $i\in\{1,2\}$, let $\varepsilon_i$ be the proportion of ordered pairs in $\Omega_i^2$ for which the equality above fails. Each triple gives six ordered pairs of distinct entries, each summing to the remaining entry. A failed pair must correspond to a triple not represented in $C$, and there are at most $r$ such triples in total. Since each $\Omega_i$ has $(2^{d/2})^2=2^d$ ordered pairs, we obtain
	\[
	(\varepsilon_1+\varepsilon_2)2^d\le6r.
	\]
	Identify each $\Omega_i$ with $\mathbb F_2^{d/2}$. By Lemma~\ref{lem:binary-linearity}, there is a linear function $\varphi_i$ on $\Omega_i$ that differs from $\psi$ at at most $\varepsilon_i2^{d/2}$ vectors.
	
	Let $\mathbf{e}_1,\ldots,\mathbf{e}_d$ be the standard basis vectors of $\mathbb F_2^d$. Define $\mathbf{u}=(u_1,\ldots,u_d)$ by
	\[
	u_j=
	\begin{cases}
		\varphi_1(\mathbf{e}_j),&1\le j\le d/2,\\
		\varphi_2(\mathbf{e}_j),&d/2<j\le d.
	\end{cases}
	\]
	By linearity, $\mathbf{u}\cdot\mathbf{x}=\varphi_i(\mathbf{x})$ for every $\mathbf{x}\in\Omega_i$. Hence
	\[
	\left|\{\mathbf{x}\in\Omega_1\cup\Omega_2:
	\mathbf{u}\cdot\mathbf{x}\ne\psi(\mathbf{x})\}\right|
	\le(\varepsilon_1+\varepsilon_2)2^{d/2}
	\le\frac{6r}{2^{d/2}}.
	\]

	We now bound $|C\setminus Q_{\mathbf{u}}|$. For $(T,X)\in C$, the set $X$ consists of exactly those vectors $\mathbf{w}$ satisfying $\mathbf{w}\cdot\mathbf{x}=\psi(\mathbf{x})$ for every entry $\mathbf{x}$ of $T$. If $(T,X)\notin Q_{\mathbf{u}}$, then $\mathbf{u}\notin X$, so $T$ contains an entry $\mathbf{x}$ with $\mathbf{u}\cdot\mathbf{x}\ne\psi(\mathbf{x})$.
	
	Fix a nonzero vector $\mathbf{x}\in\Omega_1\cup\Omega_2$. The $2^{d/2}-2$ choices of a distinct nonzero vector $\mathbf{y}$ in the same subspace count each triple containing $\mathbf{x}$ twice, by choosing either of its other two entries. Thus $\mathbf{x}$ occurs in $(2^{d/2}-2)/2$ triples. It also occurs in $2^{d/2}-1$ pairs, one for each nonzero vector in the other subspace. Therefore the total number of tuples containing $\mathbf{x}$ is
	\[
	\frac{2^{d/2}-2}{2}+(2^{d/2}-1)
	<3\cdot2^{d/2-1}.
	\]
	
	Since at most $6r/2^{d/2}$ vectors satisfy $\mathbf{u}\cdot\mathbf{x}\ne\psi(\mathbf{x})$, and $C$ contains at most one vertex for each tuple, we obtain
	\[
	|C\setminus Q_{\mathbf{u}}|
	\le\frac{6r}{2^{d/2}}\cdot3\cdot2^{d/2-1}
	=9r.
	\]
	
	For each $T\in\mathcal T$, exactly one equivalence class contains $\mathbf{u}$, so
	\(
	|Q_{\mathbf{u}}|=|\mathcal T|=|C|+r.
	\)
	Consequently,
	\[
	\begin{aligned}
		|C\triangle Q_{\mathbf{u}}|
		=|Q_{\mathbf{u}}|-|C|+2|C\setminus Q_{\mathbf{u}}|
		=r+2|C\setminus Q_{\mathbf{u}}|
		\le19r.
	\end{aligned}
	\qedhere
	\]
\end{proof}

\subsection{\texorpdfstring{The graph $H(t,d)$ and its list assignment}{The graph H(t,d) and its list assignment}}\label{subsec:twofold-construction}

For the given $d$, we fix an even-sized set $B$ as in Lemma~\ref{lem:fourier-parity}.
Let
\[
\gamma(\mathbf{u})=
\begin{cases}
49,&\mathbf{u}\in B,\\
48,&\mathbf{u}\in\mathbb F_2^d\setminus B
\end{cases}
,\enspace
b_v=\frac12\sum_{\mathbf{u}\in X_v}\gamma(\mathbf{u})\ (v\in V(J)),
\enspace\text{and}\enspace
\ell=\frac12\sum_{\mathbf{u}\in\mathbb F_2^d}\gamma(\mathbf{u}).
\]
Since $|B|$ is even, $\ell$ is an integer. We have
$24\cdot2^d\le\ell\le49\cdot2^{d-1}$.
By Observation~\ref{obs:auxiliary-graph-properties}, $|V(J)|<6\cdot2^d<\ell$.

Let $t$ be a positive integer and let $H(t,d)$ be the graph obtained
from $J$ by adding a clique $A$ and joining every vertex of $A$ to every
vertex of $J$, where $A$ consists of $\ell+4t-|V(J)|$ new vertices. By the construction, we have the following observation.

\begin{observation}\label{obs:twofold-gadget}
  $|V(H(t,d))|=\ell+4t$ and $\alpha(H(t,d))=4$.
\end{observation} 

For $v\in V(H(t,d))$, let
\begin{equation}\label{eq:twofold-list-sizes}
f(v)=\begin{cases}
\lfloor b_v\rfloor,&v\in V(J),\\
\ell,&v\in A.
\end{cases}
\end{equation}

\begin{observation}\label{obs:twofold-list-bounds}
For every $v\in V(J)$,  $|V(J)|<f(v)\le b_v\le\ell$.
Also, $\sum_{v\in V(J)}b_v=\ell|\mathcal T|$.
\end{observation}

\begin{proof}
For $v\in V(J)$, $|X_v|=2^{d-2}$ by Observation~\ref{obs:dot-product-partition}.
Since $\gamma(\mathbf{u})\ge48$, $b_v\ge\frac12\cdot48|X_v|=6\cdot2^d$.
As $6\cdot2^d$ is an integer, $|V(J)|<6\cdot2^d\le\lfloor b_v\rfloor=f(v)\le b_v$.
Since $X_v\subseteq\mathbb F_2^d$ and $\gamma(\mathbf{u})>0$, $b_v\le\ell$.

For fixed $T\in\mathcal T$, the equivalence classes $X$ for $T$ partition
$\mathbb F_2^d$, so
\[
\sum_{X:\,(T,X)\in V(J)} b_{(T,X)}
=\frac12\sum_{X:\,(T,X)\in V(J)}\sum_{\mathbf{u}\in X}\gamma(\mathbf{u})
=\frac12\sum_{\mathbf{u}\in\mathbb F_2^d}\gamma(\mathbf{u})
=\ell.
\]
Summing over $T\in\mathcal T$, we obtain
$\sum_{v\in V(J)}b_v=\ell|\mathcal T|$.
\end{proof}

\begin{lemma}\label{lem:twofold-bad-lists}
There is a list assignment $L^*$ of $H(t,d)$ with $|L^*(v)|=2f(v)+t$
for every $v\in V(H(t,d))$ such that $H(t,d)$ has no $2$-fold
$L^*$-coloring.
\end{lemma}

\begin{proof}
For each $\mathbf{u}\in\mathbb F_2^d$, we add a set of $\gamma(\mathbf{u})$
colors to every list on $Q_{\mathbf{u}}\cup A$, with pairwise disjoint
sets for distinct $\mathbf{u}$. Each $v\in A$ lies in every
$Q_{\mathbf{u}}\cup A$, so its list has
$\sum_{\mathbf{u}\in\mathbb F_2^d}\gamma(\mathbf{u})=2\ell$
colors. Each $v\in V(J)$ lies in $Q_{\mathbf{u}}$ exactly when $\mathbf{u}\in X_v$,
so its list has $\sum_{\mathbf{u}\in X_v}\gamma(\mathbf{u})=2b_v$ colors.
We keep $2f(v)$ colors in each list and add the same $t$ new colors
to every list to obtain $L^*$.

A $2$-fold $L^*$-coloring would use two colors at each vertex, for a
total of $2|V(H(t,d))|$ color uses. Each color in the set for
$\mathbf{u}$ appears only in lists on the clique $Q_{\mathbf{u}}\cup A$
and can be used at most once. Each of the $t$ new colors can be used
at most four times because $\alpha(H(t,d))=4$ by
Observation~\ref{obs:twofold-gadget}. Thus the total number of color
uses is at most $2\ell+4t<2(\ell+4t)=2|V(H(t,d))|$, a contradiction.
\end{proof}

\begin{lemma}\label{lem:twofold-seed-choosable}
If $28000t<2^{d/2-9}$, then $H(t,d)$ is $f$-choosable.
\end{lemma}

\begin{proof}
Let $L$ be an $f$-list assignment of $H(t,d)$, and suppose that
$H(t,d)$ has no $L$-coloring. We construct a partial $L$-coloring greedily.
At each step, we choose any pair of uncolored nonadjacent vertices in $J$
whose lists share an unused color and assign one such color to both
vertices. We stop when no two uncolored nonadjacent vertices of $J$
share an unused color. Let $p$ be the number of pairs colored. For each
uncolored vertex $v$ of $H(t,d)$, let $L'(v)$ consist of the unused colors
in $L(v)$.

If we could choose a color from each remaining list $L'(v)$ so that
all chosen colors were distinct, we could extend the partial coloring
to an $L$-coloring of $H(t,d)$. Thus, by Hall's theorem, there is a
nonempty set $S$ of uncolored vertices with $|Y|<|S|$, where
$Y=\bigcup_{v\in S}L'(v)$.

\begin{claim}\label{cl:twofold-hall-set-sizes}
$p+1\le4t$ and $|V(J)\setminus S|\le4t+p-1$.
\end{claim}

\begin{proof}
If $S\subseteq V(J)$, then $|S|\le |V(J)|-2p$ because the $2p$ colored
vertices all lie in $J$. For any $v\in S$, we would have
$|Y|\ge|L'(v)|\ge f(v)-p>|V(J)|-2p\ge|S|$ by
Observation~\ref{obs:twofold-list-bounds}, a contradiction.
Thus $S\cap A\ne\varnothing$. For $v\in S\cap A$, we have
$|Y|\ge|L'(v)|\ge f(v)-p=\ell-p$. By Observation~\ref{obs:twofold-gadget},
there are $\ell+4t-2p$ uncolored vertices in $H(t,d)$, so
\[
\ell-p\le |Y|<|S|\le\ell+4t-2p.
\]
It follows that $p<4t$ and $|S|>\ell-p$. Since $p$ and $|S|$ are integers,
$p+1\le4t$ and $|S|\ge\ell-p+1$. Hence
\[
|V(J)\setminus S|\le|V(H(t,d))\setminus S|=\ell+4t-|S|\le4t+p-1.\qedhere
\]
\end{proof}

For each $\alpha\in Y$, let
$C_\alpha=\{v\in S\cap V(J):\alpha\in L'(v)\}$.
The set $C_\alpha$ is a clique in $J$, possibly empty. Otherwise we
could still apply the greedy rule at the start of the proof to two
nonadjacent vertices of $C_\alpha$ using $\alpha$.

\begin{claim}\label{cl:twofold-deficit-bound}
$(\ell+4t)|\mathcal T|-\sum_{\alpha\in Y}|C_\alpha|<104\cdot2^d t$.
\end{claim}

\begin{proof}
By counting the pairs $(v,\alpha)$ with $v\in S\cap V(J)$ and
$\alpha\in L'(v)$ in two ways, we have
\begin{equation}\label{eq:twofold-clique-incidences}
\sum_{\alpha\in Y}|C_\alpha|=\sum_{v\in S\cap V(J)}|L'(v)|.
\end{equation}

By Claim~\ref{cl:twofold-hall-set-sizes}, $|V(J)\setminus S|\le4t+p-1<8t$.
For each $v\in V(J)\setminus S$, we have
$b_v\le(49/2)|X_v|=49\cdot2^{d-3}<8\cdot2^d$, since $\gamma(\mathbf{u})\le49$
for all $\mathbf{u}$ and $|X_v|=2^{d-2}$ by
Observation~\ref{obs:dot-product-partition}. Hence
\begin{equation}\label{ineq:twofold-outside-S}
\sum_{v\in V(J)\setminus S}b_v<64\cdot2^d t.
\end{equation}

For $v\in S\cap V(J)$, the list $L(v)$ has size $f(v)=\lfloor b_v\rfloor$.
Since at most $p$ colors were deleted and $0\le b_v-\lfloor b_v\rfloor<1$,
\begin{equation}\label{ineq:twofold-remaining-list-deficit}
0\le b_v-|L'(v)|=(b_v-\lfloor b_v\rfloor)+(f(v)-|L'(v)|)<p+1.
\end{equation}
\begin{samepage}
	
By Claim~\ref{cl:twofold-hall-set-sizes}, $p+1\le4t$. By
Observations~\ref{obs:twofold-list-bounds}, \ref{obs:auxiliary-graph-properties},
and~\ref{obs:tuple-count}, we have $\sum_{v\in V(J)}b_v=\ell|\mathcal T|$,
$|V(J)|=4|\mathcal T|$, and $|\mathcal T|<2\cdot2^d$, respectively.
By equation~\eqref{eq:twofold-clique-incidences} and
inequalities~\eqref{ineq:twofold-outside-S} and~\eqref{ineq:twofold-remaining-list-deficit}, 
\begin{align*}
(\ell+4t)|\mathcal T|-\sum_{\alpha\in Y}|C_\alpha|
 &=4t|\mathcal T|+\sum_{v\in V(J)}b_v-\sum_{v\in S\cap V(J)}|L'(v)|\\*
 &=4t|\mathcal T|+\sum_{v\in V(J)\setminus S}b_v
   +\sum_{v\in S\cap V(J)}(b_v-|L'(v)|)\\*
 &\le4t|\mathcal T|+64\cdot2^d t+(p+1)|V(J)|\\*
 &\le4t|\mathcal T|+64\cdot2^d t+4t\cdot4|\mathcal T|\\*
 &=20t|\mathcal T|+64\cdot2^d t<104\cdot2^d t.\qedhere
\end{align*}
\end{samepage}
\end{proof}

For each $\alpha\in Y$, let $\mathbf{z}_\alpha\in\mathbb F_2^d$ be a fixed
vector such that
$|C_\alpha\triangle Q_{\mathbf{z}_\alpha}|\le19(|\mathcal T|-|C_\alpha|)$.
By Lemma~\ref{lem:clique-approximation}, such a vector exists.
Since $|Y|<|S|\le|V(H(t,d))|=\ell+4t$, by Claim~\ref{cl:twofold-deficit-bound},
\[
\sum_{\alpha\in Y}|C_\alpha\triangle Q_{\mathbf{z}_\alpha}|
\le19\left(|Y||\mathcal T|-\sum_{\alpha\in Y}|C_\alpha|\right)
<19\cdot104\cdot2^d t.
\]

For $\mathbf{u}\in\mathbb F_2^d$, let
\[
g(\mathbf{u})=2|\{\alpha\in Y:\mathbf{z}_\alpha=\mathbf{u}\}|-\gamma(\mathbf{u}).
\]
Since $\gamma(\mathbf{u})=49$ for $\mathbf{u}\in B$ and $48$ otherwise,
the value $g(\mathbf{u})$ is odd if and only if $\mathbf{u}\in B$.
For each $v\in V(J)$, since $v\in Q_{\mathbf{z}_\alpha}$ exactly when
$\mathbf{z}_\alpha\in X_v$,
\[
\frac12\sum_{\mathbf{u}\in X_v}g(\mathbf{u})
=|\{\alpha\in Y:v\in Q_{\mathbf{z}_\alpha}\}|-b_v.
\]
For $v\in S\cap V(J)$, exactly $|L'(v)|$ colors $\alpha\in Y$ satisfy
$v\in C_\alpha$, and $|L'(v)|\le f(v)\le b_v$. For $v\in V(J)\setminus S$,
no $C_\alpha$ contains $v$. Thus at most $b_v$ of the cliques $C_\alpha$
contain $v$ for each $v\in V(J)$.
Replacing $C_\alpha$ by $Q_{\mathbf{z}_\alpha}$ changes this count only
at vertices in $C_\alpha\triangle Q_{\mathbf{z}_\alpha}$. By
Claim~\ref{cl:twofold-deficit-bound},
$\ell|\mathcal T|-\sum_{\alpha\in Y}|C_\alpha|<104\cdot2^d t$.
Thus, by the triangle inequality,
\[
\sum_{v\in V(J)}\left|\sum_{\mathbf{u}\in X_v}g(\mathbf{u})\right|
\le2\sum_{\alpha\in Y}|C_\alpha\triangle Q_{\mathbf{z}_\alpha}|
  +2\left(\ell|\mathcal T|-\sum_{\alpha\in Y}|C_\alpha|\right)
<40\cdot104\cdot2^d t.
\]

For a pair $T=(\mathbf{y},\mathbf{z})\in\mathcal T$, the four sets $X$
with $(T,X)\in V(J)$ partition $\mathbb F_2^d$. On each $X$, both
$\mathbf{u}\cdot\mathbf{y}$ and $\mathbf{u}\cdot\mathbf{z}$ are constant,
so the sign $(-1)^{\mathbf{u}\cdot\mathbf{x}}$ is constant for each
$\mathbf{x}\in\{\mathbf{0},\mathbf{y},\mathbf{z},\mathbf{y}+\mathbf{z}\}$.
Thus for each of these four vectors,
\[
|\Delta_{\mathbf{x}}(g)|\le
\sum_{X:\,(T,X)\in V(J)}\left|\sum_{\mathbf{u}\in X}g(\mathbf{u})\right|.
\]
Every $\mathbf{x}\in\mathbb F_2^d$ is a sum of one vector from $\Omega_1$
and one from $\Omega_2$. We use each nonzero component as an entry of $T$,
and for a zero component we choose any nonzero entry from its subspace.
Then $\mathbf{x}$ is one of the four vectors for this pair. Summing over
the four vectors for every pair, we obtain
\[
\sum_{\mathbf{x}}|\Delta_{\mathbf{x}}(g)|
\le4\sum_{v\in V(J)}\left|\sum_{\mathbf{u}\in X_v}g(\mathbf{u})\right|
<160\cdot104\cdot2^d t.
\]
Since $160\cdot104<28000$ and $28000t<2^{d/2-9}$, we have
$160\cdot104\cdot2^d t<2^{3d/2-9}$.
By Lemma~\ref{lem:fourier-parity},
$\sum_{\mathbf{x}}|\Delta_{\mathbf{x}}(g)|\ge2^{3d/2-9}$, a contradiction.
\end{proof}

\subsection{\texorpdfstring{Construction of the graph $G$}{Construction of the graph G}}\label{subsec:twofold-uniform}

Let $k\ge25\cdot2^{48}$ be an integer, and let $d$ be the largest even
integer with $25\cdot2^d\le k$. Then $d\ge48$ and, by maximality of $d$,
we have $k<25\cdot2^{d+2}=100\cdot2^d$.
Let $t=2^{d/2-24}$. Since $d$ is even and $d\ge48$, $t$ is a positive
integer. Also, $t>\sqrt{k}/(10\cdot2^{24})$.

Let $H=H(t,d)$ and let $f$ be the function on $V(H)$ defined in
Subsection~\ref{subsec:twofold-construction}. By equation~\eqref{eq:twofold-list-sizes} and
Observation~\ref{obs:twofold-list-bounds}, for every $v\in V(H)$, we have
\[f(v)\le\ell\le49\cdot2^{d-1}<25\cdot2^d\le k.\]

Let $W$ be a clique with vertices $w_1,\ldots,w_k$. We construct $G$
from the disjoint union of $W$ and $N$ copies $H^{(1)},\ldots,H^{(N)}$
of $H$, where
\[
N=\frac{(2k+t)!}{2^kt!}.
\]
On each copy $H^{(j)}$, we use the function $f$ from $H$. For each
$v\in V(H^{(j)})$, we join $v$ to $w_1,\ldots,w_{k-f(v)}$.
There are no edges between distinct copies.

\begin{observation}\label{obs:twofold-choice-number}
$\ch(G)=k$.
\end{observation}

\begin{proof}
Given any $k$-list assignment of $G$, we color $W$ greedily. For each
$v\in V(H^{(j)})$, we delete from its list the colors used on its
$k-f(v)$ neighbors in $W$. At least $f(v)$ colors remain. Since
$28000<2^{15}$, we have $28000t<2^{15}t=2^{d/2-9}$. By
Lemma~\ref{lem:twofold-seed-choosable}, $H$ is $f$-choosable, so we can
color each $H^{(j)}$ from the remaining lists. Thus $\ch(G)\le k$.
Since $W$ is a clique on $k$ vertices, $\ch(G)=k$.
\end{proof}

Let $X$ be a set of $2k+t$ colors, and let $L(w)=X$ for every
$w\in V(W)$. In a $2$-fold coloring of $W$ from $X$, the vertices of
$W$ receive pairwise disjoint pairs of colors. After choosing pairs
for $w_1,\ldots,w_i$, we choose the pair for $w_{i+1}$ from the
$2k+t-2i$ remaining colors. Thus the number of $2$-fold colorings of
$W$ from $X$ is
\[
\prod_{i=0}^{k-1}\binom{2k+t-2i}{2}=\frac{(2k+t)(2k+t-1)\cdots(t+2)(t+1)}{2^k}=\frac{(2k+t)!}{2^kt!}=N.
\]
Let $\phi_1,\ldots,\phi_N$ be these $2$-fold colorings.

Let $L^*$ be the list assignment from Lemma~\ref{lem:twofold-bad-lists},
with its colors renamed to lie outside $X$. For each $v\in V(H^{(j)})$, let
\[
L(v)=L^*(v)\cup\bigcup_{i=1}^{k-f(v)}\phi_j(w_i).
\]
The pairs $\phi_j(w_i)$ are pairwise disjoint subsets of $X$, so
$|L(v)|=2f(v)+t+2(k-f(v))=2k+t$.

Suppose that $G$ has a $2$-fold $L$-coloring $\psi$. The restriction
of $\psi$ to $W$ is $\phi_j$ for some $j$. For each $v\in V(H^{(j)})$,
every color in $\bigcup_{i=1}^{k-f(v)}\phi_j(w_i)$ is used on a neighbor
of $v$, so
$\psi(v)\subseteq L^*(v)$. The restriction of $\psi$ to $H^{(j)}$ is
therefore a $2$-fold $L^*$-coloring, contradicting
Lemma~\ref{lem:twofold-bad-lists}. Hence $G$ is not
$(2k+t,2)$-choosable.
This proves Theorem~\ref{thm:unbounded-gap}.

\paragraph{Disclosure of AI use.}
We used AI tools to assist in developing the constructions and proofs
and in preparing the manuscript. We carefully checked all AI-generated
suggestions that we incorporated into the paper and take full
responsibility for its content.

\end{document}